\documentclass[twoside]{amsart}
\usepackage{amssymb}
\usepackage{amsmath}
\usepackage{amsthm}
\usepackage{diagrams}
\usepackage[ansinew]{inputenc}

\newtheorem{theorem}{Theorem}[section]
\newtheorem{proposition}[theorem]{Proposition}
\newtheorem{corollary}[theorem]{Corollary}
\newtheorem{lemma}[theorem]{Lemma}

\theoremstyle{definition}
\newtheorem{definition}[theorem]{Definition}
\newtheorem{example}[theorem]{Example}
\theoremstyle{remark}
\newtheorem{remark}[theorem]{Remark}

\newcommand{\Ol}{{\mathcal O}}
\newcommand{\f}{\varphi}
\newcommand{\ra}{\longrightarrow}

\newcommand{\V}{{\mathcal V}}
\newcommand{\W}{{\mathcal W}}

\newcommand{\pu}{\mathbb P^1}
\newcommand{\proj}{\mathbb P}
\DeclareMathOperator{\loc}{\mathrm{Locus}}
\DeclareMathOperator{\cloc}{\mathrm{ChLocus}}
\newcommand{\ratcurves}{\textrm{Ratcurves}^n(X)}
\newcommand{\ratcurvesx}{\textrm{Ratcurves}^n(X,x)}
\newcommand{\om}{\textrm{Hom}}
\newcommand{\Aut}{\textrm{Aut}}
\newcommand{\Univ}{\textrm{Univ}}
\DeclareMathOperator{\pic}{Pic}
\DeclareMathOperator{\Exc}{Exc}
\DeclareMathOperator{\rk}{\mathrm{rk}}

\newcommand{\cone}{\textrm{NE}}
\newcommand{\onespan}[1]{\langle [#1] \rangle}
\newcommand{\twospan}[2]{\langle [#1],[#2] \rangle}
\newcommand{\threespan}[3]{\langle [#1],[#2],[#3] \rangle}
\newcommand{\morespan}[2]{\langle [#1], \ldots, [#2] \rangle}
\newcommand{\E}{{\mathcal E}}
\newarrow{RDashto}{}{-}{}{-}{vee}

\begin{document}

\date{}

\title{Fano manifolds with long extremal rays}

\author{Marco Andreatta, Gianluca Occhetta}

\address{Dipartimento di Matematica
via Sommarive 14\\
I-38050 Povo (TN)}
\email{andreatt@science.unitn.it}
\email{occhetta@science.unitn.it}


\begin{abstract}
Let $X$ be a Fano manifold of pseudoindex $i_X$ whose Picard number is at least two
and let $R$ be an extremal ray of $X$ with exceptional locus $\Exc(R)$.
We prove an inequality which bounds the length of $R$ in terms of
$i_X$ and of the dimension of $\Exc(R)$ and we investigate the border cases.\\
In particular we classify Fano manifolds $X$ of pseudoindex $i_X$ obtained
blowing up  a smooth variety $Y$ along a smooth subvariety $T$ such that
$\dim T < i_X$.
\end{abstract}

\maketitle


\section{Introduction}

A smooth complex projective variety of dimension $n$ is called {\sf Fano} if its
anticanonical bundle $-K_X = \wedge^n TX$ is ample.
The {\sf index} of $X$, $r_X$, is the largest natural
number $m$ such that $-K_X=mH$ for some (ample) divisor $H$ on $X$, while the 
{\sf pseudoindex} $i_X$ is defined as the minimum anticanonical degree of rational 
curves on $X$ and it is an integral multiple of $r_X$. \\
The pseudoindex is related to the Picard number $\rho_X$ of $X$
by a conjecture which claims that $\rho_X(i_X-1) \le n$,
with equality if and only if $X \simeq (\proj^{i_X-1})^{\rho_X}$; this conjecture
appeared in \cite{BCDD} as a generalization of a similar one (with the index
in place of the pseudoindex) proposed by Mukai in 1988.\\
A first step towards the proof of this conjecture was made by
Wi\'sniewski in \cite{Wimu}, where he proved that
if $i_X > \frac{n+2}{2}$ then $\rho_X =1$. 
More recently several authors (\cite{BCDD}, \cite{op}, \cite{ACO}, \cite{Caspre}) 
dealt with this problem but the general case is still open.\\
In this paper we investigate a related problem.\\
Let $X$ be a Fano manifold with $\rho_X > 1$ and let $R$ be an extremal ray of $X$.
Let $l(R) := \min\{ -K_X \cdot C~|~ C \hbox {\ \ a rational curve in } R\}$
be the length of $R$ and $\Exc(R) =\{x \in X : x \in C \hbox {\ \ a rational curve in } R\}$ 
be its exceptional locus.  We first prove the following bound:
$$
   i_X+ l(R) \leq \dim \Exc(R) +2,\eqno{(*)}
$$
which is an improved statement of the conjecture in the case $\rho_X = 2$. 

Then we investigate the cases in which equality holds.
Equivalently we ask if on a Fano variety of pseudoindex $i_X$ an extremal ray $R$
of maximal length does determine the structure of the variety. We prove the following

\begin{theorem}\label{uguale} Let $X$ be a Fano manifold of dimension $n$, pseudoindex $i_X$ and 
Picard number $\rho_X \ge 2$, and let $R$ be a fiber type or divisorial extremal ray such that
$$i_X +l(R)= \dim \Exc(R) +2.$$
Then $X \simeq \proj^k \times \proj^{n-k}$
or $X \simeq Bl_{\proj^t}(\proj^n)$ with $0 \le t \le \frac{n-3}{2}$.
\end{theorem}

We do not know how to prove a similar theorem if $R$ is an extremal ray whose associated
contraction is small (i.e. $\dim \Exc(R) \le n-2$).
However if we replace in the assumptions the pseudoindex $i_X$ with the index $r_X$ then 
we have the following

\begin{theorem}\label{ugualeindice} Let $X$ be a Fano manifold of dimension $n$,
index $r_X$, and Picard number $\rho_X \ge 2$, and let $R$ be an extremal ray such that
$$r_X +l(R)= \dim \Exc(R) +2.$$
Then, denoted by $e$ the dimension of $\Exc(R)$, we have 
$X= \proj_{\proj^k}(\Ol^{\oplus e-k+1} \oplus \Ol(1)^{\oplus n-e})$,
with $k = n-r+1$.
\end{theorem}

Finally we consider the next step, namely the case
$$i_X +l(R)= \dim \Exc(R) +1.$$

For a fiber type or divisorial extremal ray $R$ we prove that 
$\rho_X \leq 3$, describing the Kleiman-Mori cone of $X$ and
classifying the varieties with $\rho_X = 3$, (Theorem \ref{menouno}).\\

If we assume moreover that $R$ is the ray associated to a smooth blow-up,
we have a complete classification:

\begin{theorem}\label{blowup}
Let $X$ be a Fano manifold and let $R$ an extremal ray whose associate contraction $\f_R:X \to Y$
is the blow up of a smooth subvariety $T \subset Y$, such that
$$i_X +l(R) \ge n \qquad \textrm{or~equivalently} \qquad i_X \ge \dim T+1.$$
Then $X$ is one of the following
\begin{enumerate}
\item[a)] $Bl_{\proj^t}(\proj^n)$, with $\proj^t$ a linear subspace of dimension 
$\le \frac{n}{2}-1$.
\item[b)] $Bl_{\proj^t}(\mathbb Q^n)$, with $\proj^t$ a linear subspace of dimension 
$\le \frac{n}{2}-1$.
\item[c)] $Bl_{\mathbb Q^t}(\mathbb Q^n)$, with $\mathbb Q^t$ a smooth quadric 
of dimension $\le \frac{n}{2}-1$ not contained in a linear subspace of $\mathbb Q^n$.
\item[d)] $Bl_p(V_d)$ where $V_d$ is $Bl_Y(\proj^n)$ and $Y$ is a submanifold
of dimension $n-2$ and degree $\le n$ contained in an hyperplane $H$ such that $p \not \in H$.
\item[e)] $Bl_{\pu \times \{p\}}(\pu \times \proj^{n-1})$.
\end{enumerate}
\end{theorem}

Note that if $T$ is a point the condition $i_X \ge \dim T+1=1$ is empty. 
In this case the theorem is actually the main theorem of \cite{BCW}, where
Fano varieties which are the blow-up at a point of a smooth variety
are classified (those varieties correspond to cases a) and b) with $t=0$ and d)
of the above theorem). That paper has been for us a very important source of inspiration.\\

In the appendix we propose a slight variation of a result of \cite{Bo} 
relating the pseudoindexes of two Fano manifolds one of which 
is the image of the other through a birational contraction.


\section{Background material}

In (2.1) and (2.2) we recall basic definitions and facts concerning Fano-Mori contractions
and families of rational curves; our notation is consistent with the one in
\cite{Kob} to which we refer the reader.\\
Afterwards, in (2.3), for the reader's convenience we recall some results of \cite{ACO}
and \cite{CO} which are frequently used in the rest of the paper.

\subsection{Fano-Mori contractions} Let $X$ be a smooth complex Fano variety
of dimension $n$ and let $K_X$ be its canonical divisor.
By the Cone Theorem the cone of effective 1-cycles which is contained in
the $\mathbb R$-vector space of 1-cyles modulo numerical equivalence,
$\cone(X) \subset N_1(X)$, is polyhedral; a face of $\cone(X)$ is
called an {\sf extremal face} and an extremal face of dimension one
is called an {\sf extremal ray}.\\

From the structure of the cone follows that

\begin{lemma}\label{diveff} \cite[Lemme 2.1]{BCW} Let $X$ be a Fano variety and $D$ an 
effective divisor on $X$. Then there exists an extremal ray $R \subset \cone(X)$
such that $D \cdot R >0$.
\end{lemma}

To an extremal face $\sigma$ is associated a morphism with connected fibers 
$\f_\sigma:X \to W$ onto a normal variety, which contracts the curves 
whose numerical class is in $\sigma$; $\f_\sigma$ is called an {\sf extremal contraction}
or a {\sf Fano-Mori contraction}.\\
A Cartier divisor $H$ such that $H = \f_\sigma ^*A$ for an ample divisor $A$ on $W$
is called a {\sf good supporting divisor} of the map $\f_\sigma$ (or of the face $\sigma$).\\
An extremal ray $R$ (and the associated extremal contraction $\f_R$) 
is called {\sf numerically effective} ({\sf nef} for short)
or of {\sf fiber type} if $\dim W < \dim X$, 
otherwise the ray (and the contraction) is {\sf non nef} or {\sf birational}.
This terminology is due to the fact that there exists an effective
divisor $E$ such that $E \cdot R <0$ if and only if the ray is not nef.
If the codimension of the exceptional locus of a birational ray $R$ is equal to one
the ray and the associated contraction are called {\sf divisorial},
otherwise they are called {\sf small}.


\subsection{Families of rational curves} 

Let $X$ be a normal projective variety and let $\om(\pu,X)$ be the scheme parametrizing
morphisms $f: \pu \to X$. We consider the open subscheme $\om_{bir}(\pu,X) \subset \om(\pu,X)$, 
corresponding to those morphisms which are birational onto their image, and its normalization
$\om^n_{bir}(\pu,X)$. The group $\Aut(\pu)$ acts on $\om^n_{bir}(\pu,X)$ and the
quotient exists.\par

\begin{definition}
The space $\ratcurves$ is the quotient of $\om^n_{bir}(\pu,X)$ by $\Aut(\pu)$, and the space
$\Univ(X)$ is the quotient of the product action of $\Aut(\pu)$ on
$\om^n_{bir}(\pu, X)\times\pu$.
\end{definition}

\begin{definition} \label{Rf}
We define a {\sf family of rational curves} to be an irreducible component
$V \subset \ratcurves$.
Given a rational curve $f:\pu \to X$ we will call a {\sf family of
deformations} of $f$ any irreducible component $V \subset
\ratcurves$ containing the equivalence class of $f$.
\end{definition}

Given a family $V$ of rational curves, we have the following basic diagram:
\begin{diagram}[size=2em] \label{diagram}
   p^{-1}(V)=:U & \rTo^{i} & X \\
   \dTo<p & & \\
   V & &\\
\end{diagram}
where $i$ is the map induced by the evaluation $ev:\om^n_{bir} (\pu, X)\times \pu \to X$
and $p$ is a $\pu$-bundle. We define $\loc(V)$ to be the image of $U$ in $X$;
we say that $V$ is a {\sf covering family} if $\overline{\loc(V)}=X$.\par
\smallskip
If we fix a point $x \in X$ everything can be repeated starting from the scheme
$\om(\pu,X; 0 \mapsto x)$ parametrizing morphisms $f:\pu \to X$ which send $0 \in \pu$
to $x$. Given a family $V \subseteq \ratcurves$, we can consider the subscheme 
$V \cap \ratcurvesx$ parametrizing curves in $V$ passing through $x$;
we usually denote by $V_x$ a component of this subscheme.

\begin{definition} 
Let $V$ be a family of rational curves on $X$. Then
\begin{itemize}
\item[(a)] $V$ is {\sf unsplit} if it is proper; 
\item[(b)] $V$ is {\sf locally unsplit} if for the general $x \in \loc(V)$ every
component $V_x$ of $V \cap \ratcurvesx$ is proper;
\item[(c)] $V$ is {\sf generically unsplit} if there is at most a finite number of curves of $V$
passing through two general points of $\loc(V)$.
\end{itemize}
\end{definition}

\begin{proposition}\cite[IV.2.6]{Kob}\label{iowifam} 
Let $X$ be a smooth projective variety
and $V$ a family of rational curves.
Assume either that $V$ is generically unsplit and $x$ is a general point in $\loc(V)$ or that
$V$ is unsplit and $x$ is any point in $\loc(V)$ or that $x$ is a point such that
$V_x$ is unsplit. Then
 \begin{itemize}
      \item[(a)] $\dim X + \deg V \le \dim \loc(V)+\dim \loc(V_x) +1$;
      \item[(b)] $\deg V \le \dim \loc(V_x)+1$.
   \end{itemize}
\end{proposition}

This last proposition, in case $V$ is the unsplit family of deformations of a minimal extremal
rational curve, i.e. a curve of minimal degree in an extremal face of $X$, 
gives the {\sf fiber locus inequality}:

\begin{proposition}\cite{Io}, \cite{Wi1}\label{fiberlocus} Let $\f$ be a Fano-Mori contraction
of $X$ and let $E = E(\f)$ be its exceptional locus;
let $S$ be an irreducible component of a (non trivial) fiber of $\f$. Then
$$\dim E + \dim S \geq \dim X + l -1$$
where 
$$l =  \min \{ -K_X \cdot C\ |\  C \textrm{~is a rational curve in~} S\}.$$
If $\f$ is the contraction of a ray $R$, then $l$ is called the {\sf length of the ray}.
\end{proposition}

\begin{definition}\label{CF}
We define a {\sf Chow family of rational curves} to be an irreducible component
$\V \subset \textrm{Chow}(X)$ parametrizing rational and connected 1-cycles.
If $V$ is a family of rational curves, the closure of the image of
$V$ in $\textrm{Chow}(X)$ is called the {\sf Chow family associated to} $V$.\\
 We say that $V$ is {\sf quasi-unsplit} if every component of any reducible cycle in $\V$ is
numerically proportional to $V$.
\end{definition}

Let $X$ be a smooth variety, $\V^1, \dots, \V^k$ Chow families of rational curves
on $X$ and $Y$ a subset of $X$.\par

\begin{definition}
We denote by $\loc(\V^1, \dots, \V^k)_Y$ the set of points  that can be joined
to $Y$ by a connected chain of $k$ cycles belonging \underline{respectively} to 
the families $\V^1, \dots, \V^k$.\\
We denote by $\cloc_m(\V^1, \dots, \V^k)_Y$ the set of points  that can be joined to $Y$ by a
connected chain  of at most $m$ cycles belonging to the families $\V^1, \dots, \V^k$.
\end{definition}

\begin{definition}
We define a relation of {\sf rational connectedness with respect to $\V^1, \dots, \V^k$}
on $X$ in the following way: $x$ and $y$ are in rc$(\V^1,\dots,\V^k)$ relation if there
exists a chain of rational curves in $\V^1, \dots ,\V^k$ which joins $x$ and $y$, i.e.
if $y \in \cloc_m(\V^1, \dots, \V^k)_x$ for some $m$.
\end{definition}

To the rc$(\V^1,\dots,\V^k)$ relation we can associate a fibration, at least on an open subset:

\begin{theorem}\cite{Cam},\cite[IV.4.16]{Kob} \label{rcvfibration}
There exist an open subvariety $X^0 \subset X$ and a proper morphism with connected fibers
$\pi:X^0 \to Z^0$ such that
   \begin{itemize}
      \item[(a)] the rc$(\V^1,\dots,\V^k)$ relation restricts to an equivalence relation on $X^0$;
      \item[(b)] the fibers of $\pi$ are equivalence classes for the rc$(\V^1,\dots,\V^k)$ relation;
      \item[(c)] for every $z \in Z^0$ any two points in $\pi^{-1}(z)$ can be connected by a chain
      of at most $2^{\dim X - \dim Z}-1$ cycles in $\V^1, \dots, \V^k$.
   \end{itemize}
\end{theorem}

The geometry of Fano varieties is strongly related to the properties of families of rational
curves of low degree. The following is a fundamental theorem, due to Mori:

\begin{theorem}\cite{Mo2}
Through every point of a Fano variety $X$ there exists a rational curve of anticanonical degree
$\leq \dim X +1$.
\end{theorem}

\begin{remark}\label{mdf}
The families $\{V^i \subset \ratcurves \}$ containing rational curves with degree $\leq \dim X +1$
are only a finite number, so, for at least one index $i$, we have that $\loc(V^i)=X$. Among these
families we choose one with minimal anticanonical degree, and we call it a {\sf minimal covering
family}. 
\end{remark}

Let $X$ be a Fano variety and $\pi:X^0 \to Z^0$ a proper surjective morphism on a 
smooth quasiprojective variety $Z^0$ of positive dimension.\\
A relative version of Mori's theorem, \cite[Theorem 2.1]{KoMiMo}, states that, for a 
general point $z \in Z^0$, there exists a rational curve $C$ on $X$ of anticanonical 
degree $\le \dim X +1$ which meets $\pi^{-1}(z)$ without being contained in it 
(an {\sf horizontal curve}, for short).\\
As in remark \ref{mdf} we can find a family $V$ of horizontal curves such that $\loc(V)$
dominates $Z^0$ and $\deg V$ is minimal among the families
with this property. Such a family is called a {\sf minimal horizontal dominating family}
for $\pi$.

\begin{lemma}\cite[Lemma 6.5]{ACO} \label{horizontal}
Let $X$ be a Fano variety, let $\pi:X \rDashto Z$ be the fibration associated to 
a rc$(\V^1,\dots,\V^k)$ relation and let $V$ be a minimal horizontal dominating
family for $\pi$. Then
   \begin{itemize}
      \item[(a)] curves parametrized by $V$ are numerically independent from curves contracted
      by $\pi$;
      \item[(b)] $V$ is locally unsplit;
      \item[(c)] if $x$ is a general point in $\loc(V)$ and $F$ is the fiber containing $x$, then
         $$\dim (F \cap \loc(V_x))=0.$$
   \end{itemize}
\end{lemma}

\subsection{Chains of rational curves, numerical equivalence and cones}

In this subsection we present some results concerning the dimension, the maximum number
of numerically independent curves and the cone of curves of  subsets
of the form $\loc(V^1, \dots, V^k)_Y$ or $\cloc(V^1, \dots, V^k)_Y$
when $V^1, \dots, V^k$ are unsplit families and $Y$ is chosen in a suitable way.

\begin{definition} Let $V^1, \dots, V^k$ be unsplit families on $X$.
We will say that $V^1, \dots, V^k$ are numerically independent if the numerical classes
$[V^1], \dots ,[V^k]$ are linearly independent in the vector space $N^1(X)$.
If moreover $C \subset X$ is a curve we will say that $V^1, \dots, V^k$ are numerically 
independent from $C$ if in $N^1(X)$ the class of $C$ is not contained in the
vector subspace generated by $[V^1], \dots ,[V^k]$.   
\end{definition}

The following lemma is a generalization of proposition \ref{iowifam} and of
\cite[Theoreme 5.2]{BCDD}

\begin{lemma} \cite[Lemma 5.4]{ACO} \label{locy}
Let $Y \subset X$ be a closed subset and $V$ an unsplit family.
Assume that curves contained in $Y$ are numerically independent from curves in $V$, and that
$Y \cap \loc(V) \not= \emptyset$. Then for a general $y \in Y \cap \loc(V)$
\begin{itemize}
      \item[(a)] $\dim \loc(V)_Y \ge \dim (Y \cap \loc(V)) + \dim \loc(V_y);$
      \item[(b)] $\dim \loc(V)_Y \ge \dim Y + \deg V - 1$.
\end{itemize}
Moreover, if $V^1, \dots, V^k$ are numerically independent
unsplit families such that curves contained in $Y$ are numerically independent
from curves in $V^1, \dots, V^k$ then either $\loc(V^1, \ldots, V^k)_Y=\emptyset$
or
\begin{itemize}
      \item[(c)] $\dim \loc(V^1, \ldots, V^k)_Y \ge \dim Y +\sum \deg V^i -k$.
\end{itemize}
\end{lemma}

{\bf Notation}: Let $S$ be a subset of $X$. 
We write $N_1(S)=\morespan{V^1}{V^k}$ if the numerical class in $X$ of every curve 
$C \subset S$ can be written as $[C]= \sum_i a_i [C_i]$, with $a_i \in \mathbb Q$ and $C_i \in V^i$.
We write $\cone(S)=\morespan{V^1}{V^k}$  (or $\cone(S)=\morespan{R_1}{R_k}$) if the numerical 
class in $X$ of every curve $C \subset S$ can
be written as $[C]= \sum_i a_i [C_i]$, with $a_i \in \mathbb Q_{\ge 0}$ 
and $C_i \in V^i$ (or $[C_i]$ in $R_i$).\par

\begin{lemma} \cite[Lemma 1]{op} \label{numequns}
Let $Y \subset X$ be a closed subset and $V$ an unsplit family of rational curves. 
Then every curve contained in $\loc(V)_Y$ is numerically equivalent to
a linear combination with rational coefficients
   $$\lambda C_Y + \mu C_V,$$
where $C_Y$ is a curve in $Y$, $C_V$ belongs to the family $V$ and $\lambda \ge 0$.
\end{lemma}

\begin{corollary}\cite[Corollary 4.4]{ACO}\label{rhobound}
If $X$ is rationally connected with respect to some (quasi) unsplit families 
$V^1, \dots, V^k$ then $N_1(X)=\morespan{V^1}{V^k}$.
\end{corollary}

\begin{proposition}\cite[Corollary 4.2]{ACO}, \cite[Corollary 2.23]{CO}\label{coneS}~\par
\begin{enumerate}
\item[(a)]  Let $V$ be a quasi-unsplit family of rational curves and $x$ a point in $\loc(V)$.
Then $\cone(\cloc_m(V)_x)=\onespan{V}$.
\item[(b)] 
Let $V$ be a family of rational curves and $x$ a point in $X$ such that $V_x$ is unsplit.
Then $\cone(\loc(V_x))=\onespan{V}$.
\item[(c)]  Let $\sigma$ be an extremal face of $\cone(X)$, $F$ a fiber
of the associated contraction and $V$ an unsplit family independent from $\sigma$.
Then $\cone(\cloc_m(V)_F) = \twospan{\sigma}{V}$.
\end{enumerate}
\end{proposition}

\begin{corollary} \label{onedim} Let $D \subset X$ be an effective divisor and $V$ an unsplit
family numerically independent from curves in $D$ such that $D \cdot V>0$; then, for every
$x \in \loc(V)$ we have $\dim \loc(V_x)=1$; in particular, if $V$ is the family of deformations
of a minimal extremal rational curve in a ray $R$ then every non trivial fiber
of $\f_R$ is one dimensional.  
\end{corollary}

{\bf Proof.} \quad Since $D \cdot V >0$, for every $x \in \loc(V)$ we have
$D \cap \loc(V_x) \not = \emptyset$, and so $\dim (D \cap \loc(V_x)) \ge \dim \loc(V_x)-1$.
It follows that $\dim \loc(V_x)=1$, since a curve in the intersection would be
a curve in $D$ not independent from $V$.\qed \par


\section{Some technical results}

In order to make the exposition clearer, we collect in this section
some technical lemmata we will use in the proofs of the main theorems.

\begin{lemma}\label{codimension}
Let $X$ be a smooth projective variety, $A \subset X$ a subvariety of dimension $m$ 
and $V$ a covering family of rational curves for $X$.\\
Suppose that for a point $x \in X \setminus A$ the family $V_x$ is unsplit
and that $\dim \loc(V_x) \ge m+2$. Then the general curve in $V_x$ does not meet $A$.
\end{lemma}

{\bf Proof.} \quad We can assume that $V_x$ is irreducible; consider the diagram

$$
\begin{diagram}[height=2em,w=2em]
U_x & \rTo^i & X \\
\dTo<p& & \\
V_x & & \\ 
\end{diagram}
$$\par
\smallskip
and the inverse image $i^{-1}(A)$; since $V_x$ is unsplit the map $i$ is
finite to one away from $i^{-1}(x)$, so it is finite to one when restricted to
$i^{-1}(A)$ and so $\dim i^{-1}(A) = \dim A = m$.
It follows that $i^{-1}(A)$ cannot dominate $V_x$ which has dimension 
$= \dim \loc(V_x)-1 \ge m+1$. \qed \par

\begin{corollary}\label{blowzero}
Suppose that a Fano variety $X$ is the blow up of a Fano variety $Y$ along
a smooth subvariety $T$ such that $\dim T \le i_X +l(R) -3$ and denote by $E$ 
the exceptional divisor.\\
Then, if $V_Y$ is a minimal dominating family of rational curves for $Y$
and $V^*$ is a family of deformations of the strict transform of a general
curve in $V_Y$ we have $E \cdot V^*=0$.
\end{corollary}

{\bf Proof.} \quad Suppose by contradiction that $E \cdot V^* >0$.
By the canonical bundle formula we have
$$-K_Y \cdot V_Y = -K_X \cdot V^* + l(R) E \cdot V^*,$$
hence 
$$-K_Y \cdot V_Y \ge i_X +l(R) \ge \dim T+3.$$
We can therefore apply lemma \ref{codimension}, and obtain that the general curve of 
$V_Y$ does not meet $T$, a contradiction with $E \cdot V^*>0$.\qed \par

\begin{lemma}\label{unique} Let $X$ be a Fano variety whose cone of curves
is generated by a divisorial extremal ray $R_1$  with exceptional locus $E$ and a fiber type 
extremal ray $R_2$, and let $V$ be a quasi unsplit covering family of rational curves.\\
Then $[V] \in R_2$; in particular $E \cdot V >0$.
\end{lemma}

{\bf Proof.} \quad Consider the rc${\mathcal V}$ fibration 
$X \rDashto Z$. 
By proposition \ref{coneS} we have $\dim Z>0$ since $V$ is quasi unsplit and $\rho_X = 2$.\\
Then $V$ is extremal by \cite[Lemma 2.28]{CO} since $X$ has not small contractions.\\
The last assertion follows from lemma \ref{diveff}, since $E \cdot R_1<0$.\qed\par

\begin{lemma}\label{splitV}
Let $X$ be a Fano variety of dimension $n$ and pseudoindex $i_X \ge 2$ whose cone of curves 
is generated by a divisorial extremal ray $R_1$ with exceptional locus $E$ and by a 
fiber type extremal ray $R_2$.
Suppose that $l(R_1)+i_X \ge n$ and that there exists a covering family $V$ of rational curves 
of degree $\le n+1$ such that $E \cdot V=0$.\\
Then $V$ is not quasi unsplit and all 
the reducible cycles in the associated Chow family $\V$ have two 
irreducible components, $C_1$ and $C_2$, where $C_1$ and $C_2$ are curves in the rays $R_1$
and $R_2$ respectively. 
\end{lemma}

{\bf Proof.} \quad 
First of all we note that, since $E \cdot V=0$, by lemma \ref{unique}
$V$ is not quasi unsplit.\\
Let $C =\sum C_i$ be a reducible cycle in $\V$. At least one of the components of 
$C$, let it be $C_1$, has negative 
intersection with $E$; in fact, if $E \cdot C_i = 0$ for every $i$ the effective divisor $E$ 
would be numerically trivial on the whole $\cone(X)$ since $\rho_X=2$.\\
Denote by $V^1$ a family of deformations of $C_1$; if $V^1$ is not unsplit 
then there exists a reducible cycle $\sum  C_{1j}$ in $\V^1$, 
and for at least one of the components, call it $C_{11}$, we have $E \cdot C_{11} <0$.\\
Denote by $V^{11}$ a family of deformations of $C_{11}$. If $V^{11}$ is not unsplit, 
we repeat the argument, and the procedure terminates because 
$-K_X \cdot V > -K_X \cdot V^1 > -K_X \cdot V^{11} > \dots >0$.\\
Therefore every reducible cycle $\sum C_i$ in $\V$ has an irreducible component
on which $E$ is negative and such that its family of deformations is unsplit.\\
Let $\Gamma$ be one of these components and $W$ a family of deformations of
$\Gamma$; since $E \cdot \Gamma <0$ we have $\loc(W) \subset E$. We claim
that $[W] \in R_1$.
Assume by contradiction that $W$ is independent from $R_1$.\\
Denoted by $F$ a fiber of $\f_{R_1}$ meeting $\loc(W)$, by lemma \ref{locy} we have
$$n-1 \ge \dim \loc(W)_F \ge i_X-1+\dim F \ge n-1.$$ 
This forces $\loc(W)=E$, so $F \subset \loc(W)$ and we can
apply part a) of lemma \ref{locy} and get 
$$ n-1= i_X-1 +\dim F = \dim \loc(W)_F \ge \dim F +\dim \loc(W_y)$$
which implies that $\dim \loc(W_y)=i_X-1$ so $W$ is covering, 
a contradiction.\\
Therefore $[W] \in  R_1$ and for every reducible cycle $\sum C_k$ in $\V$ we have
\begin{equation}\label{ineq} 
n+1 \ge -K_X \cdot V = -K_X \cdot \sum C_k  \ge l(R_1) +(k-1)i_X \ge n.
\end{equation}
Hence $k=2$ and every reducible cycle has two components, $C_1$, which belongs to 
$R_1$ and $C_2$.\par 
\medskip
From (\ref{ineq}) it also follows that $-K_Y \cdot V \ge n$, and this implies
that $V$ is not locally unsplit. To prove this fact we assume by contradiction 
that $V$ is locally unsplit.\\
If $-K_X \cdot V=n+1$, then $\rho_X=1$ by proposition \ref{coneS} b), while
if $-K_X \cdot V=n$ for a general $x \in X$ $D_x=\loc(V)_x$ is a divisor; this
divisor is zero on $R_1$ by corollary \ref{onedim}, since $R_1$ has fibers of dimension $\ge 2$
and, by proposition \ref{coneS} b), $\cone(D_x)=\onespan{V}$.\\
On the other hand $\f(D_x)$ is an effective, hence ample
divisor on $Y$, so it meets the center of the blow up which has positive dimension.
It follows that $D \cap E \not = \emptyset$; this, together with $D \cdot R_1 = 0$
implies that $D$ contains fibers of $\f_{R_1}$, a contradiction with $\cone(D_x)=\onespan{V}$.\par
\medskip
Since $V$ is covering, not locally unsplit and $\loc(V_1)=E$, 
the family $V_2$ of deformations of 
$C_2$ is a covering family; we have $-K_X \cdot V_2 \le i_X +1<2i_X$, so $V_2$ is an unsplit 
family and therefore, by lemma \ref{unique} its numerical class belongs to the ray $R_2$.\qed 


\section{A bound on the length}

\begin{lemma}\label{familyR} 
Let $X$ be a Fano manifold with $\rho_X \geq 2$, let $R$ be an 
extremal ray of $X$ and denote by $\Exc(R)$ its exceptional locus. 
Then there exists a family of rational curves $V$ independent from $R$ such that, for
some $x \in \Exc(R)$, $V_x$ is unsplit.\\
Moreover, if $R$ is not nef and  $W$ is a minimal covering family, then, among the families 
of deformations of irreducible components of cycles in $\W$, there is a
family $V$ as above and one of the following happens
\begin{enumerate}
\item[a)] $\Exc(R) \subset \overline{\loc(V)}.$
\item[b)] There exists a reducible cycle $C_R + \sum_{i=1}^k C_i$ in $\W$ with $[C_R] \in R$. 
\end{enumerate} 
\end{lemma}

{\bf Proof.} \quad If $R$ is a nef ray it's enough to choose $V$ as
the family of deformation of a minimal extremal rational curve in any ray $R_1 \not = R$,
so we can assume that $R$ is not nef.\\
Let $W$ be a minimal covering family for $X$. Note that, since $W$ is covering, it is 
certainly independent from $R$: in fact, since $R$ is not nef there exists an effective divisor
$H$ such that $H \cdot R <0$ so curves whose numerical class is in $R$
are contained in $H$.\\
If there exists $x \in \Exc(R)$ such that $W_x$ is unsplit then we are done, otherwise
for every $x \in \Exc(R)$ there exists in $\W$ a reducible cycle $\sum_1^m C_i$, 
with rational components, passing through $x$.\\
Denote by $T^i$ the families of deformations of the curves $C_i$; since the number of
such families is finite, for at least one index $j$ we have
$\Exc(R) \subset \overline{\loc(T^{j})}$.\\
If $T^{j}$ is independent from $R$ then let $W^1=T^{j}$,
otherwise let $C_{j}+ \sum_{i \not = j} C_i$ be a reducible cycle in $\W$ passing through a point
$x \in \Exc(R)$. Since $[W]=[C_j+ \sum_{i \not = j} C_i]$ is independent
from $R$ and every component which is proportional to $R$ is contained in $\Exc(R)$
there exists an irreducible component $C_k$ independent from
$R$ which meets $\Exc(R)$. In this case denote by $W^1$ the family of deformations of
$C_k$.\\
We have thus found a family $W^1$ which is independent from $R$ such that 
$\loc(W^1) \cap \Exc(R) \not=\emptyset$. Moreover either we can choose 
$W^1$ such that $\Exc(R) \subset \overline{\loc(W^1)}$ or there
exists a reducible cycle in $\W$ with one component belonging to $R$.
Let $x_1 \in \loc(W^1) \cap E$. If $W^1_x$ is unsplit we are done,
otherwise we repeat the argument.\\ 
Since $n+1 > \deg W > \deg W^1 > \dots >0$ the procedure terminates. \qed \par

\bigskip
{\bf Proof of inequality (*).} \quad Let $x \in \Exc(R)$ 
and $V$ be as in lemma \ref{familyR}.
Let $\f_R:X \to Y$ be the extremal contraction associated to $R$ and let $F_x$ be the fiber of $\f_R$ 
which contains $x$. The numerical class of every curve in $F_x$ is in $R$ and,
by proposition \ref{coneS}, b) the numerical class of every curve in $\loc(V_x)$
is proportional to $[V]$ so, since $V$ is independent from $R$,
we have $\dim \loc(V_x) \cap F_x=0$ . Moreover, 
by inequalities \ref{iowifam} and \ref{fiberlocus}, we
have $\dim \loc(V_x) \ge i_X -1$ and $\dim F_x \ge \dim X -\dim \Exc(R) +l(R)-1$.
Combining these inequalities we get

\begin{eqnarray}\label{proofprop}
\dim X  &\ge&  \dim \loc(V_x) + \dim F_x  \\
&\ge& i_X +\dim X -\dim \Exc(R) +l(R)-2\nonumber
\end{eqnarray}

which gives
$$ i_X + l(R) \le \dim \Exc(R) +2,$$
and the proposition is proved. \qed


\section{The border cases}


{\bf Proof of \ref{uguale}.} \quad First of all note that, since the length of a fiber type extremal
ray is $\le n+1$, equality holding if and only if $X \simeq \proj^n$, and the length of a birational
extremal ray is $\le n-1$, the assumptions of the theorem imply
$i_X \ge 2$.\\
Let $V$ the family given by lemma \ref{familyR},
let $x \in \Exc(R)$ be a point such that $V_x$ is unsplit and let $F_x$ be the fiber
of $\f_R$ containing $x$.  
If equality holds in (*), then equality holds everywhere in (\ref{proofprop});
in parti\-cu\-lar we have 
\begin{eqnarray}
\dim F_x=l(R)+\dim X-\dim \Exc(R)-1 \label{dimF}\\ 
\dim \loc(V_x)=i_X -1. \label{dimloc}
\end{eqnarray}
The last equality, together with inequality \ref{iowifam} yields that 
$\dim \loc(V)=n$, so $V$ is a covering family, and that $\deg V=i_X$,  so 
$V$ is unsplit.\par
\medskip
If $\f_R:X \to Y$ is of fiber type then we can apply \cite[Theorem 1]{op} to get that
$X \simeq \proj^{i_X-1} \times \proj^{l(R)-1}$.\par
\medskip
Suppose now that $\f_R:X \to Y$ is divisorial and call $E$ the divisor $\Exc(R)$.\\
By (\ref{dimF}) we have $\dim F_x=l(R)$; note that, since $V$ is covering and unsplit,
this equality holds for every $F$, hence we can apply  \cite[Theorem 5.1]{AO}
and we obtain that $\f_R$ is the blow up of $Y$ along a smooth subvariety $T$.\\
Let $F$ be any fiber of $\f_R$; by lemma \ref{locy} b) we have
$$\dim \loc(V)_F \ge \dim F+i_X-1 \ge n,$$
so, by proposition \ref{coneS} c) we have $NE(X)=\twospan{R}{R_V}$, where
$R_V$ is the ray spanned by the numerical class of $V$.\par
\medskip
The target $Y$ of $\f_R$ is a smooth variety with $\rho_Y=1$ covered by rational curves, 
hence a Fano variety; let $V_Y$ be a minimal dominating family of rational curves for
$Y$ and let $V^*$ be the family of deformations of the strict transform of a general curve
in $V_Y$. By corollary \ref{blowzero} we have $E \cdot V^*=0$, hence, by lemma \ref{splitV}, 
the family $V^*$ is not quasi unsplit and all the reducible cycles in the associated Chow 
family $\V^*$ have two irreducible components, $C_{R}$ and $C_{V}$,
where $C_R$ and $C_{V}$ are curves in the rays $R$ and $R_V$ respectively.
In particular
\begin{equation}\label{conto}
n+1 \ge -K_Y \cdot V_Y=-K_X \cdot V^*= -K_X \cdot (C_{R}+C_{V}) \ge l(R) +i_X=n+1,
\end{equation}
and $Y \simeq \proj^n$ by the proof of 
\cite[Theorem 1.1]{Ke}. (Note that the assumptions of the quoted result are different, but
the proof actually works in our case, since for a very general $y$ the pointed family 
$V_{Y_y}$ has the properties 1-3 in \cite[Theorem 2.1]{Ke}).\par
\medskip
By equation \ref{conto} we also have $-K_X \cdot C_R= l(R)$ and 
$-K_X \cdot C_{V} = i_X$, so $C_R$ and $C_{V}$ are minimal extremal rational curves;
in particular $E \cdot C_R=-1$ and therefore,
since $E \cdot V^*=0$  we have $E \cdot C_{V}=1$.\\
Let $\psi:X \to Z$ be the contraction of the ray $R_V$; we know that $E \cdot C_{V}>0$, so 
every fiber of $\psi$ meets a fiber $F$ of $\f_R$ and therefore its dimension
is $n-\dim F=i_X-1$, since fibers of different extremal ray contractions can meet 
only in points.\\
Let now $G$ be a general fiber of $\psi$; $G$ is smooth, and, by adjunction 
$$K_G+(\dim G +1)E_G=\Ol_G,$$ 
so $G$ is a projective space and $E \cap G$ is an 
hyperplane which dominates $T$. Therefore $T$ is a projective
space by \cite[Theorem 4.1]{Laz}.\par
\smallskip
The bound on the dimension of $T$ follows from the fact that $\dim T =n -l(R)-1 =i_X-2$
and $2i_X \le l(R)+ i_X = n+1$.\qed \par

\begin{theorem}\label{menouno} Let $X$ be a Fano manifold of Picard number $\rho_X \ge 2$,
and let $R$ a fiber type or divisorial extremal ray such that
$$i_X +l(R)= \dim \Exc(R) +1.$$
Then $\rho_X \le 3$ and $\rho_X=3$ if and only if $X$ is 
\begin{enumerate}
\item[a)] $\pu \times \pu \times \proj^{n-2}$.
\item[b)] $Bl_{\pu \times \{p\}}(\pu \times \proj^{n-1})$.
\item[c)] $Bl_p(V_d)$ where $V_d$ is $Bl_Y(\proj^n)$ and $Y$ is a submanifold
of dimension $n-2$ and degree $\le n$ contained in an hyperplane $H$ such that $p \not \in H$.
\end{enumerate}
If $\rho_X=2$, except for the cases
\begin{enumerate} 
\item[d)] $Bl_p(\mathbb Q^n)$.
\item[e)] $Bl_{\proj^{l(R)-2}}(\proj^n)$. 
\end{enumerate}
the cone of curves $\cone(X)$ is generated by $R$ and by a fiber type extremal ray
and moreover $i_X \ge 2$.
\end{theorem}

{\bf Proof.} \quad If $i_X= 1$ and $R$ is divisorial we have
$l(R) \ge n-1$ so, by \cite[Theorem 1.1]{AO}, $X$ is the blow
up at a point of a variety $X'$; by \cite[Theorem 1.1]{BCW} we are in case
c) or in case d).\\
If $i_X =1$ and $R$ is of fiber type then $l(R)=n$; in particular $\f_R:X \to B$ 
is equidimensional with $n-1$-dimensional fibers over a smooth curve $B$. 
The general fiber of $\f_R$ is a projective space
by \cite[Corollary 0.4]{CMS} or \cite[Theorem 1.1]{Ke}.\\
Over an open Zariski subset $U$ of $B$ the morphism $p$
is a projective bundle. 
By taking the closure in $X$ of a
hyperplane section of $p$ defined over the open
set $U$ we get a global relative hyperplane section divisor (we use
$\rho(X/B)=1$) hence $p$ is a projective bundle globally by \cite[Lemma 2.12]{Fu1}.\\
Since $X$ is a Fano manifold $B \simeq \proj^1$. Write $X=\proj_{\pu}(\oplus \Ol(a_i))$
with $0 \le a_0 \le a_i \le a_{n-1}$. A straightforward computation
shows that $X$ is Fano if and only if either all the $a_i$ are zero or all the $a_i$
but the last are zero and $a_{n-1}=1$. In the first case $X= \pu \times \proj^{n-1}$
and $i_X=2$, in the second case $X=Bl_{\proj^{n-2}}(\proj^n)$.\par
\medskip
From now on we can assume $i_X \ge 2$.\\
Let $V$ the family given by lemma \ref{familyR},
let $x \in \Exc(R)$ be a point such that $V_x$ is unsplit and let $F_x$ be the fiber
of $\f_R$ containing $x$. First of all we prove that $V$ is an unsplit family. 
In fact, if $V$ were not unsplit then $-K_X \cdot V \ge 2i_X$ and $\dim \loc(V_x)\ge 2i_X-1$.\\
In this case we would have 
$$\dim \loc(V_x)+\dim F_x \ge 2i_X-1 +n +l(R)-\dim \Exc(R) -1\ge$$
$$\ge n +i_X-1 >n$$ 
and so $\dim \loc(V_x) \cap F_x \ge 1$, a contradiction, since $V$ is independent from $R$.\par
\medskip
Now we divide the proof in two cases, according to the type of $R$.\par
\medskip
{\bf Case 1:} $R$ is nef.\par
\medskip
Recall that, according to the proof of lemma \ref{familyR}, in this
case $V$ is the family of deformations of a minimal extremal rational curve
in a ray $R_1$ different from $R$.\par
\medskip
Suppose that $R_1$ is not nef; by inequality \ref{fiberlocus}, if $F$
is a fiber of the associated contraction we have
$\dim F \ge i_X$ and, by lemma \ref{locy}
$$\dim \loc(R)_{F} \ge \dim F +l(R)-1 \ge i_X +l(R)-1 = n.$$
It follows that $\dim F=i_X$ and $X=\loc(R)_F$, so $\cone(X)=\twospan{R}{R_1}$ 
by proposition \ref{coneS} c).\\
Since $\dim F=i_X=l(R_1)$ for every fiber of the contracion associated to $R_1$,
this contraction is a smooth blow up by \cite[Theorem 5.1]{AO}.\\
We can repeat the second part of the proof of theorem \ref{uguale},
exchanging $R_1$ and $R$ and obtain that $X=Bl_{\proj^{l(R)-2}}(\proj^n)$,
so we are in case e).\par
\medskip
Suppose now that $R_1$ is nef and consider the rc$(R,R_1)$ fibration $\pi_{R,R_1}:X -->Z$.\\
Let $F$ be a general fiber of $\pi_{R,R_1}$ and $x \in F$ a point; $F$ contains
$\loc(R,R_1)_x$ which has dimension $\ge i_X+l(R)-2=n-1$ by lemma \ref{locy},
so $\dim Z \le 1$.\par
\medskip
Suppose that $\dim Z=1$ and let $V'$ be a minimal horizontal dominating family for 
$\pi_{R,R_1}$; by lemma \ref{horizontal} c) $\dim \loc(V')_x=1$ and so $-K_X \cdot V'=2 = i_X$.\\
In particular $V'$ is unsplit and, by \ref{iowifam}, covering.
We can apply \cite[Theorem 1]{op} to conclude that $X=\pu \times \pu \times \proj^{n-2}$
and we are in case a).\par
\medskip
If $\dim Z=0$ then $X$ is rc$(R,R_1)$-connected and $\rho_X=2$ by corollary \ref{rhobound}; 
in this case we clearly have $\cone(X)=\twospan{R}{R_1}$.\par
\medskip
{\bf Case 2:} $R$ is not nef.\par
\medskip

Let $W$ be a minimal covering family for $X$ and let $V$ be a family as in lemma
\ref{familyR}, chosen among the families of deformations of irreducible components
of cycles in $\W$.\par
\medskip
{\bf Step 1} \quad V is an unsplit covering family.\par
\medskip

Let $x \in \Exc(R)$ be a point such that $V_x$ is unsplit and let
$F_x$ be the fiber of $\f_R$ containing $x$. Since $V$ is independent from $R$, 
we have $\dim \loc(V_x) \cap F_x=0$, hence $\dim \loc(V_x) \le n -\dim F_x$, giving
rise to the following chain of inequalities:

$$ -K_X \cdot V -1 \le \dim \loc(V_x) \le n-\dim F_x \le n-l(R) \le i_X.$$
 
This implies that $-K_X \cdot V \le i_X+1$ and therefore that $V$ is unsplit,
since we are assuming that $i_X \ge 2$.\par
\medskip
Suppose that $V$ is not a covering family. Then, by inequality
\ref{iowifam}, $\dim \loc(V_x) \ge 2$ and therefore $E:=\Exc(R)$ is not contained in $\loc(V)$.
In fact, in this case, by lemma \ref{locy} a) we would have 
$\dim \loc(V)_F \ge \dim F + \dim \loc(V_x)= n$, a contradiction.\\
So we are in case b) of lemma \ref{familyR} and there exists a reducible
cycle $C_R + \sum C_i$  in $\W$ with $[C_R] \in R$. Hence we have
$$n \ge \-K_X \cdot W \ge -K_X \cdot (C_R + \sum_{i=1}^k C_i) \ge l(R)+k i_X \ge n +(k-1)i_X$$
forcing $-K_X \cdot W=n$ and $k=1$.\\
We have thus proved that in $\W$ there exists a reducible cycle $C_R +C_V$,
with $C_R$ in $R$ and $C_V$ in $V$.\\
Let $D = \loc(W_x)$ for a general $x \in X$; by proposition \ref{coneS} b)
$\cone(D)=\onespan{W}$.\\ 
By corollary \ref{onedim}, since the fibers of $\f_R$ are at least two dimensional we have
$D \cdot R=0$; by the same corollary, since $\dim \loc(V_x) \ge 2$ we have
$D \cdot V=0$. This implies also that $D \cdot W= D \cdot (C_R + C_V)= 0$.\\
By lemma \ref{diveff} there exists an extremal ray $R_1$ such that $D \cdot R_1 >0$; 
let $V^1$ be a family of deformations of a minimal curve in $R_1$.
By lemma \ref{locy} b) we have $\dim \loc(V^1)_D \ge \dim D+i_X-1 \ge n$,
hence $X= \loc(V^1)_{D}$ and $\rho_X=2$.\\
This is a contradiction, since $D$ is zero on $R$ and $V$ and so, if $\rho_X=2$
it would be zero on the entire cone.
Therefore $V$ is a covering family as claimed.\qed \par
\medskip
{\bf Step 2}\quad $\rho_X \le 3$.\par
\medskip
Let $F$ be a fiber of $\f_R$; by lemma \ref{locy} b) we have 
$$\dim \loc(V)_F \ge \dim F+i_X-1 \ge n-1$$
If $X = \loc(V)_F$ then, by proposition \ref{coneS} c) $\cone(X) = \twospan{R}{V}$ and
we are done. Note that this is always the case if $\dim F > l(R)$, so we assume from
now on that $\f_R$ is equidimensional with fibers of dimension $l(R)$, hence it is
a smooth blow up by \cite[Theorem 5.1]{AO}.\\
An irreducible component of $\loc(V)_F$ is thus a divisor $D \subset X$ such
that $\cone(D)=\twospan{R}{V}$. If $D \cdot V >0$ then 
$X=\cloc_2(V)_F$ and $\cone(X) = \twospan{R}{V}$ again by 
proposition \ref{coneS} c), so we can assume $D \cdot V=0$.\par
By lemma \ref{diveff} there exists an extremal ray $R_1$ such that $D \cdot R_1>0$.\par
\medskip
If $R_1 \not \subset \cone(D)$ then, by lemma \ref{locy} b), denoted by $V^1$
a family of deformations of a minimal extremal rational curve in $R_1$, we have
$\dim \loc(V^1)_D=n$. By lemma \ref{numequns} $N_1(X)=\threespan{R}{V}{V^1}$, so 
$\rho_X \le 3$, equality holding if and only if $R_1$ is not contained in the
vector subspace of $N_1(X)$ spanned by  $R$ and $[V]$.\par 
\medskip
If $R_1 \subset \cone(D)$ then $R_1=R$ because $D \cdot V=0$.
It follows that $\loc(R)_D=E$, so $N_1(E)=\twospan{R}{V}$.\\
If $E \cdot V>0$ then $\loc(V)_E=X$ and $N_1(X)=\twospan{R}{V}$ by lemma
\ref{numequns}, so $\rho_X = 2$.
We claim that we cannot have $E \cdot V=0$; in fact, in this case every curve of $V$ which meets 
$E$ is entirely contained in $E$, so $E=\loc(V)_F=D$ and we have $D \cdot R<0$.
Recalling that $D \cdot V=0$ we have that $D$ is not positive on $\cone(D)$, a
contradiction, since we are assuming $R_1 \subset \cone(D)$ and $D \cdot R_1 >0$.\par
\medskip
{\bf Step 3} \quad $\rho_X=2$, description of the cone.\par
\medskip
We have to prove that $\cone(X)=\twospan{R}{R_1}$ where $R_1$
is a fiber type extremal ray. By step two 
this is the case if for a fiber $F$ of $\f_R$ either we have $X=\loc(V)_F$ or
an irreducible component of $\loc(V)_F$ is a divisor $D$ such that
$D \cdot V>0$.
We can therefore assume that an irreducible component of 
$\loc(V)_F$ is a divisor $D$ such that $D \cdot V=0$; moreover we know that 
there exists an extremal ray $R_1$ of $X$
on which $D$ is positive.\par
\medskip
If $R_1 \not \subset \cone(D)$ then $\cone(X)=\twospan{R}{R_1}$ and moreover, by corollary
\ref{onedim} the contraction associated to $R_1$ has one dimensional fibers, and so
it is of fiber type, since $i_X \ge 2$.\par
\medskip
If $R_1  \subset \cone(D)$ then $R_1=R$ thus, if $V$ is not extremal, $D$ is negative 
on an extremal ray $R_2$, and so $\Exc(R_2) \subset D$, against $\cone(D)=\twospan{R}{V}$.
Therefore $V$ is extremal and $\cone(X)=\twospan{R}{V}$.\par
\medskip
{\bf Step 4} \quad $\rho_X=3$, description of the cone.\par
\medskip
By step two, if $\rho_X=3$, then $\loc(V)_F$ has dimension $n-1$;
moreover, denoted by $D$ one irreducible component of $\loc(V)_F$
we have $D \cdot V=0$ and $D \cdot R_1>0$ for a ray $R_1$ not contained in the
vector subspace of $N_1(X)$ spanned by  $R$ and $[V]$.\\
Since $\cone(D)=\twospan{R}{V}$, by corollary \ref{onedim}, every fiber
of the contraction associated to $R_1$ is one dimensional.
Combining this with $i_X \ge 2$, by inequality \ref{fiberlocus}, we have that
$V^1$ is a covering unsplit family.\\
By lemma \ref{locy}, denoting again by $F$ a fiber of $\f_R$ we have 
$ \dim \loc(V,V^1)_F =\dim \loc(V^1,V)_F =n$, so $X =\loc(V,V^1)_F=\loc(V^1,V)_F$.\\
We can write $X=\loc(V,V^1)_F=\loc(V)_{\loc(V^1)_F}$ and therefore, by lemma 
\ref{numequns} and proposition \ref{coneS} the numerical class of 
every curve in $X$ can be written as a linear combination $a[V]+b[V^1]+c[R]$ with $b,c \ge 0$.\\
On the other hand $X=\loc(V^1,V)_F= \loc(V^1)_{\loc(V)_F}$, so the numerical class of 
every curve in $X$ can be written as a linear combination $a[V]+b[V^1]+c[R]$ with $a,c \ge 0$.
By the uniqueness of the decomposition it follows that $\cone(X)=\threespan{V}{V^1}{R}$.\par
\medskip
{\bf Step 5} \quad If $\rho_X=3, i_X \ge 2$ and $R$ is not nef then 
$X \simeq Bl_{\pu \times \{p\}}(\pu \times \proj^{n-1}).$ \par
\medskip
We have thus proved that the cone of curves of $X$ is generated by $R$, which is the ray
associated to a smooth blow up $\f_R:X \to Y$,
and by other two fiber type extremal rays, call them $R_1$ and $R_2$, which both have
length two.  In particular we have $i_X=2$, so $l(R)=n-2$ and $\dim F_R=n-2$
for every fiber of $\f_R$.\\
Moreover, since $E=\Exc(R)$ is non negative on $R_1$ and $R_2$, 
by \cite[Proposition 3.4]{Wi1} $Y$ is a Fano variety.\\
The effective divisor $E$ is positive on at least one of the rays $R_i$
by lemma \ref{diveff}; let us assume that $E \cdot R_1 >0$.
Let $\sigma$ be the extremal face spanned by $R$ and $R_1$ and consider the associated
contraction $\f_\sigma$.\\
Let $x \in X$ be a point, let $\Gamma_1$ be a curve in $R_1$ through $x$ and 
let $F$ be a fiber of $\f_R$ meeting $\Gamma_1$.
The fiber of $\f_\sigma$ through $x$ contains $\loc(R_1)_{F}$, which 
has dimension $n-1$ by lemma \ref{locy}, so the target of $\f_\sigma$
is a smooth curve, which has to be rational since $X$ is Fano.
We have a commutative diagram
$$
\begin{diagram}
X & \rTo^{\f_R}& Y \\
 &\rdTo<{\f_\sigma}& \dTo>{\psi_\sigma} \\
& & \pu
\end{diagram}
$$

The general fiber $F_\sigma$ of $\f_\sigma$ is, by adjunction, a Fano variety of index $\ge 2$
which has a divisorial extremal ray of length $\dim F_\sigma-1$, so, by theorem \ref{uguale},
$F_\sigma \simeq Bl_p{\proj^{n-1}}$.\\
It follows that the general fiber of $\psi_\sigma$ is $\proj^{n-1}$. The Fano variety
$Y$ has a fiber type extremal ray $\psi_\sigma$ of length $\dim Y$ while
the other ray is of fiber type, since the associated contraction contracts
the images of curves in $R_2$. Therefore $i_Y \ge 2$.\\
We can thus apply theorem \ref{uguale} to conclude that $Y \simeq \pu \times \proj^{n-1}$.
Let $T \simeq \pu$ be the center of the blow up; we claim that $T$ is a fiber of
the projection $Y \to \proj^{n-1}$.
By contradiction, assume that this is not the case. Let $C \simeq \pu$ be a fiber
of the projection $Y \to \proj^{n-1}$ meeting $T$ and let $\widetilde C$
be the strict transform of $C$.\\
By the canonical bundle formula we have
$$-K_X \cdot \widetilde C = -K_Y \cdot C -l(R)E \cdot \widetilde C  \le 2-l(R) \le 0,$$
and so $X$ is not a Fano variety, a contradiction. \qed \par

\section{Blow ups}

{\bf Proof of \ref{blowup}.} \quad 
If $i_X+l(R)=n+1$, by theorem \ref{uguale} we have that $X=Bl_{\proj^{t}}(\proj^n)$, with
$t \le \frac{n-3}{2}$.\par
\medskip
We can thus assume that $i_X+l(R)=n$.
By theorem \ref{menouno}, if $\rho_X \ge 3$, then $X$ is either
$Bl_{\pu \times \{p\}}(\pu \times \proj^{n-1})$ or $Bl_p(V_d)$ where $V_d$ is $Bl_Y(\proj^n)$ 
and $Y$ is a submanifold of dimension $n-2$ and degree $\le n$ contained in an hyperplane
which does not contain $p$. Note that case a) of theorem \ref{menouno}
has been excluded since it is not a blow up.\par
\medskip
We can thus assume, from now on, that $\rho_X=2$; again by theorem \ref{menouno}
either $X \simeq Bl_p(\mathbb Q^n)$ or $\rho_X=2$, $i_X \ge 2$ and 
the cone of curves of $X$ is generated by $R$ and by a fiber type 
extremal ray $R_V$. (Case e) of theorem \ref{menouno} has been excluded
since in that case $R$ is a fiber type ray).\par
\medskip
The target $Y$ of $\f_R$ is a smooth variety with $\rho_Y=1$ covered by rational curves, 
hence a Fano variety; let $V_Y$ be a minimal dominating family of rational curves for
$Y$ and let $V^*$ be the family of deformations of the strict transform of a general curve
in $V_Y$. 
The center of the blow up $T$, has dimension $\le \dim Y-3$ since 
$${\rm codim~} T-1 =l(R) \ge i_X\ge 2,$$
therefore we can apply corollary \ref{blowzero} and obtain that 
$E \cdot V^* =0$.\par
\medskip
Since $E \cdot V^*=0$, by lemma \ref{splitV}, the family
$V^*$ is not quasi unsplit and all the reducible cycles in the associated Chow 
family $\V^*$ have two irreducible components, $C_1$ and $C_2$,
where $C_1$ and $C_2$ are curves in the rays $R$ and $R_V$ respectively.\\
Let $\Gamma_R$ and $\Gamma_V$ be curves in $R$ and $R_V$ respectively
with minimal anticanonical degree. Since $\f_R$ is a smooth blow up 
$E \cdot \Gamma_R=-1$, hence the numerical class of every curve in $R$
is an integral multiple of $[\Gamma_R]$; in particular we can write $[C_1]=m_1[\Gamma_R]$
with $m_1$ a positive integer.
By the canonical bundle formula
\begin{eqnarray} \label{kooc}
n+1 &\ge& -K_Y \cdot V_Y=-K_X \cdot V^*= -K_X \cdot (C_1+C_2) \ge\\
&\ge& m_1 l(R) + i_X \ge (m_1-1)l(R) + n. \nonumber
\end{eqnarray}

Recalling that $l(R) \ge i_X \ge 2$ we have $m_1=1$, i.e. $[C_1]=[\Gamma_R]$. It follows
that $E \cdot C_2=1$, so $[C_2]=[\Gamma_V]$ and $[V^*]=[\Gamma_R+\Gamma_V]$.\par
\medskip
Consider now the contraction of $R_V$, $\psi:X \to Z$ and let $F$ be any fiber of $\psi$.\\
Since $E \cdot \Gamma_V >0$ the fiber $F$ meets 
a fiber $F_R$ of $\f_R$ and therefore $\dim F \le n -\dim F_R=i_X$.\\
On the other hand, by inequality \ref{fiberlocus} 
$\dim F \ge l(R_V)-1 \ge i_X-1$, so the length of $R_V$ is either $i_X$ or $i_X+1$.
In the first case, by equation \ref{kooc} we have $-K_Y \cdot V_Y =n$,
while in the second we have $-K_Y \cdot V_Y =n+1$.\\
The contraction $\psi$ is supported by $K_X+i_X E$ in the first case and by
$K_X+(i_X+1)E$ in the second; in both cases, since for every fiber of 
$\psi$ we have $i_X-1 \le \dim F \le i_X$, the target variety $Z$ is smooth
by \cite[Theorem 4.1]{AWD}.\par
\medskip
The general fiber of $\psi$ has dimension either $i_X-1$ or $i_X$, so the dimension of
$Z$ is either $l(R)+1$ or $l(R)$. We divide the proof in two cases, accordingly. \par
\medskip
{\bf Case 1} \quad $\dim Z=l(R)+1$.\par
\medskip
In this case $\psi$ is supported by $K_X+i_XE$, its general fiber
has dimension $i_X-1$ and it is a projective space $\proj^{i_X-1}$ 
by \cite[Theorem 4.1]{AWD}, while jumping fibers, if they exist, have dimension $i_X$ and are
projective spaces $\proj^{i_X}$, again by \cite[Theorem 4.1]{AWD}.\par
\medskip
We claim that, for at least one fiber $F$ of $\psi$, we have $E \cap F=\proj^{i_X-1}$.
The claim is clearly true if either $E$ contains a fiber of dimension $i_X-1$ or,
being $E \cdot \Gamma_V=1$, if $\psi$ has a jumping fiber
($E$ cannot contain a jumping fiber $F$, otherwise, by lemma \ref{locy} a)
we will have $\dim E \ge \dim \loc(R)_F \ge i_X +l(R) \ge n$). \\
Suppose by contradiction that neither of these two possibilities happens.
The restriction of $\psi$ to $E$ is thus an equidimensional
morphism with general fiber a projective space, such that $E$
restricted to the general fiber is $\Ol_\proj(1)$, so $\psi$
makes $E$ a projective bundle over $Z$.\\
Therefore $E$, which is also a projective bundle over $T$, has two projective bundle 
structures and $\rho_E=2$ so, by \cite[Theorem 2]{OW}, 
$E$ is the projectivization of the tangent bundle of a projective space, but this is 
impossible since the two fibrations of $E$
have fibers of dimension $i_X-2$ and $l(R)$ and these two
dimensions are different, being $l(R) \ge i_X$, so the claim is proved.\par
\medskip
It follows that either $\psi$ has a jumping fiber or $E$ contains a fiber of $\psi$;
in both cases $T$, the center of the blow up, is dominated by the intersection of $E$
with this fiber, and so it is a projective space of dimension $i_X-1$ 
by \cite[Theorem 4.1]{Laz}.\par
\medskip
To finish the proof, we have to show that $Y \simeq \mathbb Q^n$, and we will
do this proving the existence of a line bundle $L_Y \in \pic(Y)$ such that
$-K_Y=nL_Y$ and applying the Kobayashi-Ochiai theorem.\\
Take a line $l$ in $T$ and denote by $Y_l$ the inverse image $\f_R^{-1}(l)$;
$Y_l$ is a projective bundle over a smooth rational curve, so a toric variety.
The restriction $\psi_{|Y_l}:Y_l \to Z$ is thus a surjective morphism from
a toric variety to a smooth variety with Picard number one, so $Z$ is a projective space 
by \cite[Theorem 1]{OW}.\\
Let $L$ be the line bundle $\psi^*\Ol_\proj(1)+E$; we have $L \cdot R= 0$ 
and therefore there exists ${L_Y} \in \pic(Y)$ such that $\f_R^*L_Y=L$.\\
Moreover, since $L \cdot V^*=1$ we have ${L_Y} \cdot V_Y =1$,
so, recalling that $-K_Y \cdot V_Y = \dim Y$ we get $-K_Y=nL_Y$ with 
and we conclude that $Y \simeq \mathbb Q^n$ by the Kobayashi-Ochiai theorem.\par

\medskip
{\bf Case 2} \quad $\dim Z=l(R)$.\par
\medskip
In this case, as noted above, every fiber of $\psi$ has dimension $i_X+1$.
The contraction $\psi$ is supported either by $K_X+i_XE$ and it is
a projective bundle or by $K_X+(i_X+1)E$ and it is a quadric bundle,
by \cite[Theorem 4.1]{AWD}.\\
Every fiber of $\f_R$ dominates $Z$ so, by \cite[Theorem 4.1]{Laz} $Z$
is a projective space.\par
Let $L$ be the line bundle $\psi^*\Ol_\proj(1)+E$; we have $L \cdot R= 0$ 
so there exists ${L_Y} \in \pic(Y)$ such that $\f_R^*L_Y=L$.\\
Moreover, since $L \cdot V^*=1$ we have ${L_Y} \cdot V_Y =1$.

\medskip
{\bf Case 2a} \quad $\psi:X \to Z$ is a projective bundle.\par
\medskip
In this case $-K_Y \cdot V_Y=n+1$, so $-K_Y=(n+1)L_Y$ and $Y$
is a projective space.
The intersection of $E$ with the general fiber of $\psi$ is thus 
a projective space and therefore the center $T$ of the blow up is a linear
space by \cite[Theorem 4.1]{Laz}.\par
\medskip
{\bf Case 2b} \quad $\psi:X \to Z$ is a quadric bundle.\par
\medskip
In this case $-K_Y \cdot V_Y=n$, so $-K_Y=nL_Y$ and $Y$
is a smooth quadric by the Kobayashi-Ochiai theorem.\\
The intersection of $E$ with the general fiber of $\psi$ is thus 
a smooth quadric, so the center $T$ of the blow up is either a linear
space or a smooth quadric by \cite{PS}.\par 
Actually the first case can be excluded by direct computation, since
the blow up of a quadric along a linear subspace is not
a quadric bundle over $\proj^r$.\\
In the second case let $\Pi \simeq \proj^i_X$ be the linear subspace of
dimension $i_X$ of $\proj^{n+1}$
which contains $T \simeq \mathbb Q^{i_X-1}$.\\
Two cases are possible: either $Y \supseteq \Pi$ or $Y \cap \Pi = T$.
The first case has to be excluded because, if $Y \supseteq \Pi$ the blow up of $\mathbb Q^n$ 
along $T$ does not give rise to a Fano variety.\\
To see this, take a line $l \subset \Pi$ not contained in $T$;
by the canonical bundle formula, if $X=Bl_{T}\mathbb Q^n$ we have
$$-K_X \cdot \tilde l = -K_Y \cdot l - l(R)E \cdot \tilde l \le n-2l(R) \le 0.$$\par
\medskip
Finally note that in both cases the bound on the dimension of the center follows from the fact
that $i_X \le l(R)$ and so $2 i_X \le l(R)+i_X \le n$. \qed


\section{Varieties with a polarization}


{\bf Proof of \ref{ugualeindice}.} \quad  Let $V$ the family given by lemma \ref{familyR},
let $x \in \Exc(R)$ be a point such that $V_x$ is unsplit and let $F_x$ be the fiber
of $\f_R$ containing $x$.\\ 
First of all we prove that $\rho_X=2$ and that the cone of curves of $X$ is generated
by $R$ and by the ray spanned by $[V]$.\\ 
We are assuming that equality holds in (*), so equality holds everywhere in (\ref{proofprop});
in parti\-cu\-lar we have 
\begin{eqnarray}
\dim F_x =l(R)+\dim X-\dim \Exc(R)-1 =\dim X-r_X+1\\ 
\dim \loc(V_x)=r_X -1. 
\end{eqnarray}
This forces $\deg V=r_X$, so the family
$V$ is unsplit. Moreover, by inequality \ref{iowifam} $V$ is a covering family.\\
Therefore, by lemma \ref{locy} we have 
$\dim \loc(V)_{F_x} \ge \dim F_{x}+r_X-1 = \dim X$, so, by
proposition \ref{coneS} c), we have $NE(X)=\twospan{V}{R}$.\par 
\medskip
Let $\psi:X \to Z$ be the contraction of the ray $R_V$ spanned by $[V]$, which
is of fiber type since $V$ is a covering family;
curves parametrized by $V$ have anticanonical degree $r_X$, so they are minimal extremal curves
in $R_V$ which has length $r_X$.\\
By inequality \ref{fiberlocus}, every fiber of $\psi$ has dimension $\ge l(R_V)-1 =r_X-1$,
so $\dim Z \le n-r_X+1$. Again by inequality \ref{fiberlocus} the fibers of $\f_R$ 
have dimension $\ge n-e+l-1 =n-r_X+1$, so they dominate $Z$. 
In particular every fiber of $\psi$ meets a fiber $F_R$ of $\f_R$ and so
its dimension is $ \le \dim X-\dim F_R=  r_X-1$; therefore
the contraction $\psi:X \to Z$ is equidimensional.\\ 
Moreover we also have that the dimension of every fiber of $\f_R$
is $\le \dim Z \le n-r_X+1$, so $\f_R$ is equidimensional
with fibers of dimension $n-r_X+1$ and $\dim Z=n-r_X+1$.\par
\medskip
Denote by $H$ the divisor such that $-K_X =r_X H$.
The general fiber $G$ of $\psi$ is, by generic
smoothness and adjunction, a projective space $\proj^{r_X}-1$ and $H_{G} \simeq \Ol(1)$,
so, by \cite[Lemma 2.12]{Fu1}, $\psi$ is a projective bundle over $Z$, $X=\proj_Z(\E)$,
with $\E=\f_R^*H$.
In particular $Z$ is a smooth Fano variety of Picard number one.\\
The canonical bundle formula yields
$$\psi^*(K_Z +\det \E) = K_X +r_XH =\Ol_X,$$
and so $-K_Z= \det \E$. Note also that, if $C_R$ is a curve in $R$ then 
\begin{equation}\label{degH}
H \cdot C_R = \dfrac{-K_X \cdot C_R}{r_X} \ge \dfrac{l(R)}{r_X}.
\end{equation}
Let $V_Z$ be a minimal covering family for $Z$ and $C$ a curve in $V_Z$;
Let $\nu: \pu \to C \subset Z$ be the normalization of $C$ and let
$Z_C$ be the fiber product $Z_C= \pu \times_C X$.

$$
\begin{diagram}[loose,height=2em,l>=3em]  
Z_C & \rTo^{\bar \nu} & X\\
\dTo<{p} & &\dTo>{\psi}\\
\pu & \rTo_{\nu} & Z 
\end{diagram}
$$

The variety $Z_C$ is a projective bundle over $\pu$, 
$Z_C=\proj_\pu(\nu^*\E)$; the vector bundle $\nu^*\E$ is ample,
so we can write $\nu^*\E \simeq \oplus_0^{r_X-1} \Ol(a_i)$ with $a_i >0$ and 
$a_i \le a_{i+1} \quad \forall i$.\\
Denote by $m$ the maximum index $i$ such that $a_i =a_0$ and rewrite $\nu^*\E$ in the following way
$$\nu^*\E \simeq \oplus^{m+1} \Ol(a_0) \oplus_{i=m+1}^{r_X-1} \Ol(a_i).$$
The cone of curves $\cone(Z_C)$ is generated by the class of a line in a fiber of $p$ 
and by the class of a section $C_0$ corresponding to a surjection 
$\nu^*\to \Ol(a_0)$.\\
The cone of curves $\cone(X)$ is generated by the class $[V]$ of a line in a fiber
of $\f_V$ and by the class of $\Gamma_R$, a minimal extremal curve in $R$.\\
The morphism $\bar \nu$ induces a map of spaces of cycles $N_1(Z_C) \to N_1(X)$
which allows us to identify $\cone(Z_C)$ with a subcone of $\cone(X)$.\\
Since $\bar \nu(Z_C)$ contains lines in the fibers of $\psi$ and contains
curves in the fibers of $\f_R$ (since for dimensional reasons $\dim (\bar \nu(Z_C) \cap F_R) \ge 1$),
we have an identification $\cone(Z_C) \simeq \cone(X)$.\\
In particular $F_R \cap \bar \nu(Z_C)$, which is a curve
whose numerical class in $X$ is a multiple of $[\Gamma_R]$, 
is the image of a curve $\Gamma$ whose
numerical class in $Z_C$ is a multiple of $[C_0]$.\\
By \cite[Lemma 3]{op} the curve $\Gamma$ is the union of disjoint minimal sections, 
so $\bar \nu(Z_C) \cap F_R$ consists of the images via $\bar \nu$ of disjoint minimal sections.\\
On the other hand, if $C_0$ is a minimal section, then $\bar \nu(C_0)$ is a curve
whose numerical class is in $R$, so it is contained in a fiber of $\f_R$.\\
It follows that the dimension of $\f_R(\Exc(R))$ is the dimension of the
space parametrizing minimal sections, which is $m$.
Therefore
$$ m= \dim \Exc(R) - \dim F_R  = l(R) +r_X -2 -\dim F_R.$$
Moreover, since $[C_0] \in R$ we have, by equation \ref{degH}
$$a_0 = H \cdot C_0 \ge \dfrac{l(R)}{r_X},$$
hence $a_i \ge a_0+1\ge \dfrac{l(R)}{r_X}+1$ for $i=m+1,\dots ,r_X-1$.\\
It follows that
$$\dim Z+1 \ge-K_Z \cdot C= \det \E\cdot C=$$
$$=(m+1)a_0 + \sum_{m+1}^{r_X-1} a_i \ge   (m+1) \dfrac{l(R)}{r_X} +(r_X-m-1)
\left(\dfrac{l(R)}{r_X}+1 \right)=$$
$$ = l(R) +r_X -m-1 = \dim F_R+1 =\dim Z+1.$$
Therefore $Z$ admits a minimal dominating family of degree $\dim Z+1$, hence 
$Z$ is a projective space of dimension $n-r_X+1$ by the proof of \cite[Theorem 1.1]{Ke}.\\
Since equality holds everywhere we also have $a_0=1$, $a_i=2 \quad i=m+1, \dots, r_X-1$,
so the splitting type of $\E$ on lines of $Z$ is uniform.\\
If $\dim \Exc(R) \le  \dim X-2$ then $\rk \E =r_X \le l(R) < n-r+1 = \dim Z$, therefore
$\E$ is decomposable by \cite{EHS} and $\E \simeq \oplus^{m+1} \Ol(1)\oplus^{r_X-1-m} \Ol(2)$.\\
If $\dim \Exc(R) = \dim X-1$ then $\rk \E =\dim Z$ and the splitting type
of $\E$ is $(1, \dots,1,2)$, so, by \cite{EHS}, either $\E$ is decomposable
or $\E$ is the tangent bundle of $Z=\proj^{\dim Z}$, but the second case
has to be excluded since $X$ has a divisorial contraction.\\
Finally, if $\Exc(R) = X$ then the splitting type
of $\E$ is $(1, \dots,1)$, so $\E$ is decomposable by \cite[Proposition 1.2]{AWI} 
and $X$ is a product of projective spaces. \qed\par

\begin{proposition}  
Let $X$ be a Fano variety of Picard number $\rho_X = 2$,
index $r_X \ge 2$, and let $R$ a fiber type or divisorial extremal ray such that
$r_X +l(R)= \dim \Exc(R) +1$. Then, if $R$ is divisorial either
$X$ is as in theorem \ref{blowup} or $X$ has the structure of a 
projective bundle over a smooth variety.\\
If $R$ is of fiber type then $X$ is a 
projective bundle or a quadric bundle or  the projectivization of a B{\v a}nic{\v a} sheaf
over a smooth variety $Y$.
\end{proposition}

{\bf Proof.} \quad By theorem \ref{menouno} either 
$Bl_{\proj^{l(R)-2}}(\proj^n)$ or the cone of curves $\cone(X)$
is generated by $R$ and by a fiber type extremal ray; let $\psi:X \to Z$ be the contraction 
of this ray.\\
Let $H$ be the line bundle such that $-K_X =r_XH$, let $A \in \pic(Z)$ be an ample divisor 
and let $H'=H+\psi^*A$. The contraction $\psi$ is supported by $K_X +r_XH'$.\par
\medskip
If $R$ is divisorial then every fiber of $\f_R$ has dimension $\ge l(R)$.
If equality holds for every fiber, $\f_R$ is a smooth blow up by  \cite[Theorem 5.1]{AO},
so $X$ is as in theorem \ref{blowup}.\\
We can therefore assume that there exists a fiber $F$ of $\f_R$ of dimension $\ge l(R)+1$.\\
The contraction $\psi:X \to Z$ has fibers of dimension $\ge r_X -1 \ge n- l(R)-1$,
so $\dim Z \le l(R)+1$. It follows that $F$ dominates $Z$ and meets every fiber of $\psi$,
forcing the equidimensionality of $\psi$.\\
We can now conclude that $X$ is a projective bundle over $Z$ by \cite[Lemma 2.12]{Fu1} since
$H \cdot V=1$.\par
\medskip
If $R$ is of fiber type then every fiber of $\f_R$ has dimension $\ge l(R)-1$
and so the contraction $\psi:X \to Z$ has fibers of dimension $\le n- l(R)+1 \le r_X$,
so we can conclude by \cite[Theorem 4.1]{AWD} and \cite[Proposition 2.5]{BW}.\qed \par

\section{Appendix}

The results in theorem \ref{blowup} show that if a Fano variety $X$ is the blow-up
of a smooth variety $Y$ along a smooth subvariety $T$ and $i_X \ge \dim T +1$
then also $Y$ is a Fano variety and $i_Y \ge i_X$.\\
In general these two facts are not true; in \cite[Section 3]{Wi1} the question
whether $Y$ has to be a Fano variety was posed and some answers were given in
\cite[Propositions 3.4 and 3.6]{Wi1}.\\
In particular the examples \cite[3.7, 3.8]{Wi1} show that $i_T \ge \dim T+1$ is the best
possible bound which guarantees that $Y$ is a Fano manifold.\par
\medskip
The second problem, i.e. - assuming that $Y$ is Fano can the pseudoindex of $Y$
be less than the pseudoindex of $X$? - has been studied in \cite{Bo}.
The following example of that paper shows that the answer can be positive:\par

\begin{example} Let $Y_n = \proj_{\proj^m}(\Ol^{\oplus m} \oplus 
\Ol(1))$ and let $T_n \subset Y_n$ be the submanifold defined by the 
subbundle  $\Ol^{\oplus m}$. Note that $\dim Y_n := n = 2m$ and $\dim T_n = m$.
Let $\f_n: X_n = Bl_{T_n}(Y_n) \ra Y_n$ be the blow-up of $Y_n$ along $T_n$.\\
One can easily prove that  $X_n$ and $Y_n$ are Fano manifolds, if $n \geq 4$,
and moreover $i_{Y_n} = 1$ and, if $n \geq 6$ ,
$i_{X_n} = 2$ (while $i_{X_4} = 1$).
\end{example}

The following are the main results of \cite{Bo}:

\begin{theorem} Let $X$ be a Fano manifold of dimension $n$
which is the blow up $\f_R :X \ra Y$ of a smooth Fano manifold $Y$ 
along a smooth subvariety $T$.
\begin{itemize}
\item If $2 \dim T \le n +i_Y -1$ then $i_X \le i_Y$ unless $n \ge 6$ is even and
$X=X_n$, $Y=Y_n$, $T=T_n$ are as in the above example.
\item If $i_Y \ge \frac{n}{3}-1$ then $i_X \le i_Y$ unless $n = 6$  and
$X=X_6$, $Y=Y_6$, $T=T_6$ are as in the above example.
\end{itemize}
\end{theorem}

We propose here a slight variation of the results of \cite{Bo},
considering birational contractions between smooth Fano manifolds:

\begin{proposition}\label{boge} Let $X$ be a Fano manifold, let $\f_R :X \ra Y$ be the
contraction of a birational extremal ray $R$ such that $Y$ is a smooth Fano manifold
and let $T = \f_R(\Exc(R))$.
If $i_Y > 2 \dim T +1 -n$ or if $i_Y > \frac{n}{3} -1$
then $i_X \leq i_Y$.
\end{proposition}

{\bf Proof.}\quad  Since $\f$ is a birational map between smooth varieties
the exceptional locus $\Exc(R)$ is a divisor and we have the formula:
$$K_X = K_Y + G, $$
where $G$ is a divisor supported on $\Exc(R)$

Let $C \subset Y$ be a rational curve such that   $i_Y = -K_Y \cdot C$
and let $V$ be a family of rational curves on $Y$ containing $C$. 

By inequality \ref{iowifam} we have 
$2 \dim \loc(V) \geq n+i_Y -1$, therefore  if $i_Y > 2 \dim T +1 -n$
we have that $\dim \loc(V)  > \dim T$ and this implies that there 
exists a curve $C$ in $V$ not contained in $T$.\\
The strict transform of it, call it $\widetilde C$, is a rational curve
on $X$ satisfying $G \cdot \widetilde C \geq 0$, therefore, by the canonical bundle formula,
$i_X \leq -{K_X} \cdot \widetilde C \leq i_Y$.\par
\medskip
Assume now that $i_Y > \frac{n}{3} -1$
and, by contradiction, that $i_X  > i_Y$;
by the first part we can assume that $i_Y \leq  2 \dim T +1 -n$.\\
Denote by $F$ a general fiber of the map $\f$; 
from  \ref{iowifam} we have $\dim F \ge i_X$ and therefore
$$i_Y \leq  2 \dim T +1 -n = n - 1 - 2 \dim F  \leq n-1-2i_X \leq n-3 -2i_Y$$
that is $i_Y \leq  \frac{n}{3} -1$, which is a contradiction.\qed \par

\begin{corollary} Let $X$ be a Fano manifold, let $\f_R :X \ra Y$ be the
contraction of a birational extremal ray $R$ such that $Y$ is a smooth Fano manifold
and let $T = \f_R(\Exc(R))$.
If $i_X \ge \dim T$ then $i_X \leq i_Y$.
\end{corollary}

{\bf Proof.} \quad Let $F$ be a general fiber of $\f_R$; we have 
$$\dim T \le \dim \Exc(R)- \dim F \le n-1 -i_X \le n -\dim T +1$$
so that 
$$2 \dim T-n +1\le 0 < i_Y$$
We can thus apply proposition \ref{boge} to conclude. \qed \par

\end{document}